\documentclass[a4paper, 10 pt, conference, twocolumn]{ieeetran} 
\IEEEoverridecommandlockouts
 \pdfoutput=1
\usepackage[utf8]{inputenc}
\usepackage[T1]{fontenc}
\usepackage{url}
\usepackage{amsmath}
\usepackage{amsfonts}
\usepackage{amssymb}
\usepackage{algorithm}
\usepackage{algpseudocode}
\usepackage{bbm}
\usepackage{graphicx}
\usepackage{caption}
\usepackage{subcaption}
\usepackage[usenames, dvipsnames]{color}
\usepackage{xcolor}
\usepackage{ifthen}
\usepackage{mathtools}
\usepackage{comment}
\newtheorem{theorem}{Theorem}%[section]

\newtheorem{assumption}{Assumption} %[section]

\newtheorem{remark}{Remark}

\newtheorem{example}{Example}
\usepackage{flushend} %to get a nice last page
\usepackage{pgfplots}

\usepackage[left=1.4cm,right=1.4cm,top=2.4cm,bottom=3.7cm]{geometry}

\newboolean{showcomments}
\setboolean{showcomments}{true}

\newcommand{\david}[1]{\ifthenelse{\boolean{showcomments}}
{\textcolor{Blue}{(David says: #1)}}{}}
\newcommand{\emma}[1]{\ifthenelse{\boolean{showcomments}}
{\textcolor{VioletRed}{(Emma says: #1)}}{}}
\newcommand{\edit}[1]{\ifthenelse{\boolean{showcomments}}
{\textcolor{Black}{#1}}{}}

\title{\LARGE \bf Optimal Control of Linear Cost Networks}

\author{ {David Ohlin, Emma Tegling and Anders Rantzer} 
\thanks{The authors are with the Department of Automatic Control and ELLIIT Strategic Research Area at Lund University, Lund, Sweden. Email: \{{\tt\small{david.ohlin, anders.rantzer, emma.tegling}\}@control.lth.se}}\thanks{This work is partially funded by Wallenberg AI, Autonomous Systems and Software Program (WASP) and the European Research Council (ERC) under the European Union’s Horizon 2020 research and innovation programme under grant agreement No 834142 (ScalableControl).}}

\begin{document}
\maketitle

\begin{abstract}
    We present a method for optimal control with respect to a linear cost function for positive linear systems with coupled input constraints. We show that the optimal cost function and resulting sparse state feedback for these systems can be computed by linear programming. Our framework admits a range of network routing problems with underlying linear dynamics. These dynamics can be used to model traditional graph-theoretical problems like shortest path as a special case, but can also capture more complex behaviors. We provide an asynchronous and distributed value iteration algorithm for obtaining the optimal cost function and control law.
\end{abstract}

\section{Introduction}

In many real-world systems the quantities of interest, like amounts of goods or distributions of probabilities, are intrinsically positive. Such systems naturally lead to models that are positive systems, where the state $x$ is confined to the positive orthant~$\mathbb{R}^n_+$. In the context of optimal control, positive systems exhibit many advantageous properties, as reviewed in~\cite{rantzervalcher18tutorial}. Of special note is that such systems admit linear Lyapunov functions. This guarantees positive stage costs for a linear cost function, giving significant computational advantages in applications where such a cost function can accurately capture the desired objective. In contrast, current methods for optimal control of general linear systems with quadratic cost~\cite{bertsekas18stable} give rise to Riccati equations where the number of variables scales quadratically with the dimension of the state, as opposed to linearly. Additionally, as noted in \cite{rantzervalcher18tutorial}, systems on this form allow for sparse optimal state feedback $u = Kx$, giving favorable scaling results for high-dimensional inputs.

%In this paper we exploit these advantages for the purpose of optimal control of a system where the magnitude of combinations of input signals are jointly constrained by a linear function of the state. These coupled constraints give rise to trade-offs between several different inputs, where in the previous work \cite{rantzer22explicit} all inputs are fully actuated in either the positive or negative direction. This novel possibility covers a range of problems excluded by the requirements in \cite{rantzer22explicit}. In particular, it gives a natural framework for a general class of network routing and shortest path problems with underlying linear dynamics unfolding on a graph. This setup, like in the quadratic cost case, allows for an explicit solution to the associated Bellman equation. Using linear costs for the state and input, a method for obtaining a linear optimal cost function $J^*(x)$ by means of linear programming is provided. 

In this paper we exploit these advantages for the purpose of optimal control of a system where the magnitude of combinations of input signals are jointly constrained by a linear function of the state. \edit{These coupled constraints generalize the previous work~\cite{rantzer22explicit}, wherein} all inputs are fully actuated in either the positive or negative direction, \edit{and give rise to trade-offs between several different inputs.} This novel \edit{formulation} covers a range of problems excluded by the \edit{constraints} in~\cite{rantzer22explicit}. In particular, it \edit{provides} a natural framework for a general class of network routing and shortest path-like problems with underlying linear dynamics unfolding on a graph. \edit{Our formulation, which uses linear costs for the state and input, allows} for an explicit solution to the associated Bellman equation \edit{just like in the quadratic cost case. We also provide} a method for obtaining the optimal cost function $J^*(x)$ by means of linear programming. 

This paper focuses on control problems over networks, in which physical constraints and scalability concerns in applications typically call for distributed solution methods. We give a distributed implementation of the value iteration algorithm for the proposed problem formulation, similar to those treated in \cite{bertsekas20multiagent}, and show that convergence results presented in \cite{bertsekas22abstract} can be applied to the present setting. This allows for distributed and asynchronous computation of the local optimal state feedback, wherein nodes need only share their local estimates of the cost function and the cost of the optimal control action. An extended analysis of the problem setup in \cite{rantzer22explicit}, leveraging established theory of Dynamic Programming and extending the problem to constraints based on other norms, can be found in \cite{li23exact} wherein upper bounds on the convergence rate of policy iteration for such systems are given. This analysis, however, also fails to capture the class of network routing problems treated here.

The formulation we present connects results from the field of optimal control with a rich history of works (e.g. \cite{dijkstra59note},~\cite{hart68heuristic}) on graphs with finite state and action spaces, which becomes a special case of the problem as stated here. Our framework is, however, flexible enough to cover a wider range of problems with continuous state and action spaces. Notably, interpreting the continuous state as a probability distribution allows for the modeling of Markov Decision Processes (MDPs) as another special case where existing results on convergence of methods like policy iteration \cite{ye11simplex} can be applied. Most works in the related literature that treat dynamic versions of the shortest path problem use models where the structure of the graph changes \cite{gao12dynamic}. Similar problems with dynamics in a continuous state space evolving on a graph structure are, however, apart from the special case of MDPs, relatively unexplored.

In the next section we present relevant existing theory and define notation. Section III contains the problem statement, required assumptions and main result in the form of Theorem~\ref{th:main}. In Section IV, an algorithm is given for distributed and asynchronous value iteration that finds the fixed point of the operator given in Theorem \ref{th:main}, corresponding to the optimal cost function. The results are illustrated by an example in Section V, where the results are applied to a simplified linear model of a large-scale cooling system and the corresponding optimal routing. Finally, Section VI presents conclusions and sets out possible avenues for future study.

\newgeometry{left=2cm,right=1.4cm,top=2.4cm,bottom=3.7cm}

\section{Preliminaries}

\subsection{Notation} 
Inequalities are applied element-wise for matrices and vectors throughout. Further, the notation $\mathbb{R}^n_+$ is used to denote the closed nonnegative orthant of dimension $n$. The operator~$\min\{A,0\}$ extracts the minimum element of $A$, yielding zero if $A$ has no negative elements. Let $\textnormal{diag}(a)$ denote a square matrix with the elements of the vector $a$ along the diagonal. The expressions $\mathbf{1}_{p\times q}$ and $\mathbf{0}_{p\times q}$ signify a matrix of ones or zeros, respectively, of the indicated dimension, with subscript omitted when the size is clear from the context. If the dimension is zero, this is to be interpreted as the empty matrix.

\subsection{Problem setup}

Consider the infinite-horizon optimal control problem
\begin{equation}
    \label{eq:optprob}
    \begin{aligned}
        \textnormal{Minimize} &\;\;\; \sum\limits_{t=0}^{\infty}\left[ s^\top x(t) + r^\top u(t) \right] \;\textnormal{over}\; \{u(t)\}^\infty_{t=0}\\
        \textnormal{subject to} &\;\;\; x(t+1) = Ax(t) + Bu(t)\\
        &\;\;\; u(t) \ge 0, \;\;\; x(0) = x_0\\
        & \;\;\; \begin{matrix} \mathbf{1}^\top u_1(t) & \le & E_1^\top x(t) \\ \vdots & & \vdots \\ \mathbf{1}^\top u_M(t) & \le & E_M^\top x(t) \end{matrix}
    \end{aligned}
\end{equation}
where $A\!\in \mathbb{R}^{n\times n}$ and $B\!=\!\begin{bmatrix}B_1 ~\cdots~ B_M\end{bmatrix}\!\in\!\mathbb{R}^{n\times m}$ with $B_i\in\mathbb{R}^{n\times m_i}$ define the linear dynamics. The input signal~$u\in\mathbb{R}^n$ is partitioned into~$M$ subvectors $u_i$, each containing $m_i$ elements, so that $m = \sum_{i = 0}^{M} m_i$. In the upcoming treatment, we will let the matrices $B_i$ be further subdivided into individual columns as $B_i\!=\!\begin{bmatrix} B_{i1}~\cdots~B_{im_i} \end{bmatrix}$ with~$B_{ij}\in\mathbb{R}^n$. The costs connected to the states and actions are $s \in \mathbb{R}^{n}_{\ge 0}$ and~$r \in \mathbb{R}^{m}_{\ge 0}$ with $r_i\in\mathbb{R}^{m_i}$ following the partition of $u$, and $r_{ij}$ denoting the $j$th element of $r_i$. The constraints on the input signal~$u(t)$ are given by $E\!=\!\begin{bmatrix}E_1 ~\cdots~ E_M\end{bmatrix}^\top\!\in\!\mathbb{R}^{M\times n}_+$. 

Solving the problem (\ref{eq:optprob}) is equivalent to finding some optimal value function $J^*(x)$ that satisfies the corresponding Bellman equation
\begin{equation}\label{eq:bellman}
    J^*(x) = \min_{u\in U(x}\left[g(x,u) + J^*(Ax + Bu)\right].
\end{equation}
In the specific context of (\ref{eq:optprob}) the immediate cost is given by the value at each instant $t$ of the cost function to be minimized, so $g(x(t),u(t)) = s^\top x(t) + r^\top u(t)$. In the previous work \cite{rantzer22explicit} the Bellman equation (\ref{eq:bellman}) has a unique linear solution $J^*(x) = p^\top x$ with $p\in\mathbb{R}^n$ for a certain class of positive linear time-invariant systems. Our main result below proves that similarly favorable results can be obtained for systems on the form (\ref{eq:optprob}), which have a natural interpretation as network routing problems but are not in general admissible in \cite{rantzer22explicit}.

\section{Main result}
In anticipation of our main result, we introduce two assumptions on the dynamics of the system in (\ref{eq:optprob}).
\begin{assumption}
\label{as:B}
        Any negative elements of $B_i$ are on its $i$th row.
\end{assumption}
% Assume that $B$ admits a decomposition $B = B^+ + B^-$ where $B^+\in \mathbb{R}^{n\times m}_+$ and 
%     \begin{equation*}
%         B^- = \begin{bmatrix} B^-_1 & & \\ & \ddots & \\ & & B^-_M \end{bmatrix}
%     \end{equation*}
%     with $B^-_i\in \mathbb{R}^{1\times m_i}$ and all other elements of $B^-$ are zero.
\begin{remark}
    Assumption \ref{as:B} is required to guarantee positivity of the dynamics in (\ref{eq:optprob}). While it may seem restricting at first glance, the partition of $u$ to allow for such a decomposition is natural in many relevant applications. For example, in the case of a network problem where the inputs $u$ correspond to the quantities transferred over each edge in the network, Assumption \ref{as:B} translates into the requirement that all edges leading from a single node $i$ must be included in the same constraint partition $u_i$. Since positivity of the system is often motivated by the states corresponding to real quantities, this only means that no more than the quantity that is already in node $i$ at time $t$ (or is expected to be within one time step) can be transferred away from it.
\end{remark}
\begin{assumption}
    \label{as:A}
    The matrices $A$, $B$ and $E$ satisfy
    \begin{equation*}
        A \ge -\textnormal{diag}\left(\begin{bmatrix} \min\{B_1,0\} \cdots \min\{B_M,0\} \end{bmatrix}\right)E \ge 0.
    \end{equation*}
\end{assumption}
\vspace{2mm}

It should be noted that these assumptions are necessary and sufficient to ensure invariance of $\mathbb{R}^n_+$ in the general case, but in special cases where parts of the state space are unreachable they may be relaxed. Given these assumptions on the matrices of the problem (\ref{eq:optprob}), we are now ready to state our main result:
\begin{theorem}\label{th:main}
    Under Assumptions \ref{as:B} and \ref{as:A}, the following three statements are equivalent:
    \begin{itemize}
        \item[($i$)] The problem (\ref{eq:optprob}) has a finite value for every $x_0 \in \mathbb{R}^{n}_+$.        
        \item[($ii$)] There exists $p \in \mathbb{R}^{n}_+$ satisfying the equation
        \begin{equation}
        \label{eq:p}
            p = s + A^\top p + \sum_{i = 1}^M \min \{ r_i +B_i^\top p, 0 \} E_i.
        \end{equation}
        \item[($iii$)] The value of the linear program
        \begin{align*}
            \textnormal{Maximize} &\;\;\; \mathbf{1}^\top p \;\textnormal{over}\; p \in \mathbb{R}^{n}_+\\
        \textnormal{subject to} &\;\;\; p \le s + A^\top p - \sum_{i = 1}^M z_i E_i\\
        &\;\;\; z_i \ge r_{ij} + B_{ij}^\top p \;\;\; \textnormal{for} \;\;\; j = 1,...,m_i\\
        &\;\;\; z_i \ge 0
        \end{align*}
        is bounded.
    \end{itemize}
    Given the existence of a $p$ as in ($ii$), the minimal value of~\eqref{eq:optprob} is given by $p^\top x_0$. The optimal linear state feedback law is then given by $u_i(t) = K_ix(t)$ with
    \begin{equation}
        K_i := \begin{bmatrix} \mathbf{0}_{j-1\times n} \\ E_i^\top \\ \mathbf{0}_{m_i-j\times n} \end{bmatrix} \;\textnormal{for}\;i=1,...,M,
    \end{equation}
    where the vector $E_i^\top$ enters at the $j$th row with $j$ being the index of the minimal element of $r_i+B_i^\top p$, provided it is negative. If all elements are nonnegative then $K_i = \mathbf{0}_{m_i\times n}$.
    % \begin{equation*}
    %     K_i := \begin{bmatrix} \begin{cases}
    %         E_i^\top & \left[(r + B^\top p)_i\right]_1 = \min\{(r + B^\top p)_i, 0\}\\
    %         0 & \text{otherwise}
    %     \end{cases} \\ \vdots \\ \begin{cases}
    %         E_i^\top & \left[(r + B^\top p)_i\right]_{m_i} = \min\{(r + B^\top p)_i, 0\}\\
    %         0 & \text{otherwise}
    %     \end{cases} \end{bmatrix}.
    % \end{equation*}
\end{theorem}
\vspace{2mm}
\begin{proof}
    We first show that $\mathbb{R}^n_+$ is invariant under the dynamics in (\ref{eq:optprob}). It holds that
    \begin{equation*}
        Ax + Bu = Ax + \begin{bmatrix} B_1\!&\!\cdots\!&\!B_n \end{bmatrix}\begin{bmatrix} u_1 \\ \vdots \\ u_M \end{bmatrix} = Ax + \sum_{i=1}^M B_iu_i.
    \end{equation*}
    As a consequence of Assumption \ref{as:B} and the constraints given in \eqref{eq:optprob} this expression has a lower bound given by
    \begin{multline*}
        Ax + \sum_{i=1}^M B_iu_i \ge Ax + \begin{bmatrix} \min\{B_1,0\} E_1^\top \\ \vdots \\ \min\{B_M,0\} E_M^\top \end{bmatrix}x\\
        = (A + \textnormal{diag}\left(\begin{bmatrix} \min\{B_1,0\} \cdots \min\{B_M,0\} \end{bmatrix}\right)E)x
    \end{multline*}
    recalling that $\min\{B_i,0\}$ denotes the minimum element of~$B_i$, yielding zero if all elements are nonnegative. This final expression is contained in $\mathbb{R}^n_+$ due to Assumption \ref{as:A}.
    %where $B_i\in \mathbb{R}^{1\times m_i}$ is a vector of the elements in the $i$:th row of $B$ that multiply $u_i$. The inequality above holds as a consequence of Assumption \ref{as:B}. Further
    %\begin{align}
        %&Ax + \begin{bmatrix} B_1 u_1 \\ \vdots \\ B_M u_M \end{bmatrix} \ge Ax + \begin{bmatrix} \min\{B_1,0\} E_1^\top x \\ \vdots \\ \min\{B_M,0\} E_M^\top x \end{bmatrix} \label{eq:positive} \\ &= (A + \textnormal{diag}\left(\begin{bmatrix} \min\{B_1,0\} & \cdots & \min\{B_M,0\} \end{bmatrix}\right)E)x,\nonumber
    %\end{align}
    %which is contained in $\mathbb{R}^n_+$ due to Assumption \ref{as:A}. 
    
    Next, let $J_k(x) = p_k^\top x$. Given the ansatz $p_0 = 0$, we compute each subsequent iterate $p_{k+1}$ by applying the Bellman equation
    \begin{align*}
        &\min_{u\in U(x)} \left[g(x,u) + J_k(f(x,u))\right]\\ 
        &= \min_{u\in U(x)} \left[s^\top x + r^\top u + p_k^\top(Ax + Bu)\right]\\
        &= s^\top x + p_k^\top Ax + \min_{u\in U(x)} \left[(r + B^\top p_k)^\top u\right]
    \end{align*}
    where $U(x)$ denotes the set of all inputs $u$ satisfying the constraints in (\ref{eq:optprob}). We partition the vector $r +B^\top p_k$ in the same fashion as $u$ into $M$ subvectors $r_i+B_i^\top p_k\in \mathbb{R}^{m_i}$. Since each subvector $u_i$ is only bound by one of the upper input constraints, we can minimize the expression above independently for each $u_i$:
    \begin{align}
    \label{eq:sum}
        \min_{u\in U(x)} (r + B^\top p_k)^\top u &= \sum_{i = 1}^M \min_{\substack{\mathbf{1}^\top u_i \le E_i x \\ u \ge 0}} (r_i+B_i^\top p_k)^\top u_i\\ \nonumber
        &= \sum_{i = 1}^M \min \{ r_i+B_i^\top p_k, 0 \} E_i^\top x \nonumber
    \end{align}
    where the minimum in the final expression is taken over all elements in the vector $r_i + B_i^\top p_k$ and 0. The Bellman iteration becomes
    \begin{align*}
        &\min_{u\in U(x)} \left[g(x,u) + J^*(f(x,u))\right]\\
        &= s^\top x + p_k^\top Ax + \sum_{i = 1}^M \min \{ r_i+B_i^\top p_k, 0 \} E_i^\top x\\
        &= \left(s + A^\top p_k + \sum_{i = 1}^M \min \{ (r_i+B_i^\top p_k, 0 \} E_i\right)^\top x\\
        &= p_{k+1}^\top x.
    \end{align*}
    As a consequence, we have an increasing sequence of functions generated by iteration of the Bellman equation, and the cost functions all take the form $J_k(x) = p_k^\top x$. Given that the problem (\ref{eq:optprob}) has a finite value according to ($i$), the sequence of value functions generated by this iteration has the upper limit $J^*(x) = p^\top x$. This limit fulfills the equation in ($ii$). The converse also holds since, given that a $p$ satisfying ($ii$) exists, the increasing sequence of $\{p_k\}_{k=0}^{\infty}$ must satisfy $p_k \le p$ for all $k$. This means that the Bellman equation has a nonnegative and finite solution, implying that ($i$) holds. To find the optimal controller we seek the policy $\mu(x)$ minimizing the Bellman equation for the optimal cost function $J^*(x) = p^\top x$. Taking the expression from (\ref{eq:sum}) we get
    \begin{equation*}
        \mu(x) = \textnormal{arg}\min_{u\in U(x)} (r + B^\top p)^\top u
    \end{equation*}
    which, partitioned like the minimization in (\ref{eq:sum}), can be performed separately for each subvector $\mu_i$ yielding the expression
    \begin{align*}
        \mu_i(x) = \underbrace{\begin{bmatrix} \mathbf{0}_{j-1\times n} \\ E_i^\top \\ \mathbf{0}_{m_i-j\times n} \end{bmatrix}}_{=K_i}x.
    \end{align*}
    where the nonzero row vector $E_i^\top$ enters on the $j$th row and $j\in\{1,...,m_i\}$ is the index of the largest in magnitude negative element of~$r_{i}~+~B_{i}^\top p$. If all elements are positive then $\mu_i(x) = \mathbf{0}_{m_i\times n}$ is the minimizing argument. As a result, the state feedback law for each partition $u_i$ is stationary and given by $u_i(t) = K_ix(t)$.
    
    Finally, we address ($iii$) in the same way as treated in \cite{rantzer22explicit}. If~($ii$) holds then the $p$ fulfilling Equation \ref{eq:p} clearly also solves the program
    \begin{align*}
            \textnormal{Maximize} &\;\;\; \mathbf{1}^\top p \;\textnormal{over}\; p \in \mathbb{R}^{n}_+\\
        \textnormal{subject to} &\;\;\; p \le s + A^\top p + \sum_{i = 1}^M \min\{r_i+B_i^\top p,0\} E_i
    \end{align*}
    which is equivalent to the linear program in ($iii$). Any solution fulfilling the Bellman inequality gives a lower bound on the corresponding optimal cost function. Thus the solution of the linear program must also be the optimal cost function and, as a consequence, ($ii$) implies ($iii$). To show the converse, assume that the linear program in ($iii$) is bounded by some argument $p_k$. Calculate
    \begin{equation*}
        p_{k+1} = s + A^\top p_k + \sum_{i = 1}^M \min \{ r_i+B_i^\top p_k, 0 \} E_i.
    \end{equation*}
    We now have $p_{k+1}\ge p_k$ and as a consequence
    \begin{align*}
        p_{k+1} &= s + A^\top p_{k} + \sum_{i = 1}^M \min \{ r_i+B_i^\top p_k, 0 \} E_i\\
        &= \min_{u\in U(x)} \left[s^\top x + r^\top u + p_{k}^\top(Ax + Bu)\right]\\
        &\le \min_{u\in U(x)} \left[s^\top x + r^\top u + p_{k+1}^\top(Ax + Bu)\right]\\
        &= s + A^\top p_{k+1} + \sum_{i = 1}^M \min \{ r_i+B_i^\top p_{k+1}, 0 \} E_i
    \end{align*}
    which is exactly the constraint in the reformulated program above. Since we assumed that $p_k$ attained the maximum of the program and $p_{k+1} \ge p_k$, the only possible conclusion is that~$p_{k+1} = p_k$. Thus the maximizing argument of the program in ($iii$), if it is bounded, also solves Equation \ref{eq:p}, implying~($ii$). This concludes the proof.
\end{proof}

To illustrate the above theorem we apply it to a simple shortest-path problem. A detailed example showing more general possibilities is presented in Section V. 

\begin{figure}[h!]
    \centering
    \includegraphics[width=.6\linewidth]{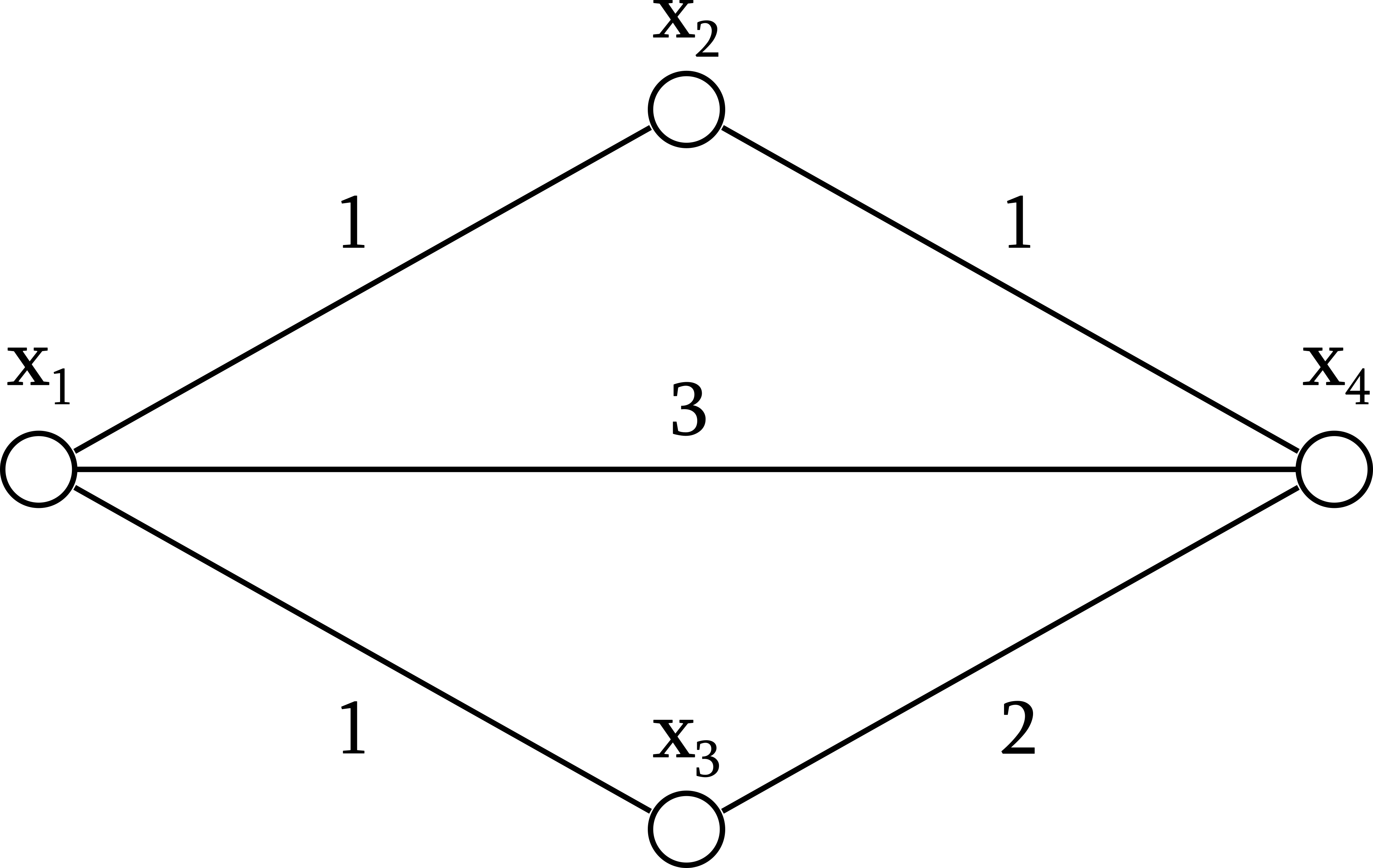}
    \caption{Weighted graph in Example \ref{ex:toy}.}
    \label{fig:toynet}
\end{figure}

\begin{example}\label{ex:toy}
    Consider the (trivial) problem of finding the shortest path from node $x_1$ to node $x_4$ in Figure \ref{fig:toynet}. This can be formulated as an optimal control problem on the form~(\ref{eq:optprob}) with~$n = 4$, $m = 10$ (one input for each direction of traversal of the edges) in the following way: Let the states~$x_i$ correspond to the nodes, select $E = A = I$ and let $B$ be a node-link incidence matrix:
    \begin{equation*}
        B = \footnotesize{\begin{bmatrix} -1 & -1 & -1 & 0 & 1 & 0 & 1 & 1 & 0 & 0 \\ 
                            1 & 0 & 0 & -1 & -1 & 0 & 0 & 0 & 1 & 0 \\
                            0 & 1 & 0 & 0 & 0 & -1 & -1 & 0 & 0 & 1 \\
                            0 & 0 & 1 & 1 & 0 & 1 & 0 & -1 & -1 & -1 \\ 
        \end{bmatrix}}.
    \end{equation*}
    We partition the input vector $u = \left[u_1^\top\;u_2^\top\;u_3^\top\;u_4^\top\right]^\top$ to collect all actions corresponding to directed edges originating in node~$i$ into $u_i$. The weights are selected as
    \begin{equation}\label{eq:weights}
        s^\top\!=\!\begin{bmatrix} 1\!&\!1\!&\!1\!&\!0 \end{bmatrix}, \; r^\top\!=\!\begin{bmatrix} 0\!&\!0\!&\!2\!&\!0\!&\!0\!&\!1\!&\!0\!&\!2\!&\!0\!&\!1 \end{bmatrix}.
    \end{equation} 
    We now have the problem on the form required by Theorem~\ref{th:main} and finding a $p$ satisfying (\ref{eq:p}) by fixed point iteration we get the optimal cost function as
    \begin{equation*}
        J^*(x) = \begin{bmatrix} 2 & 1 & 2 & 0 \end{bmatrix}^\top x.
    \end{equation*}
    The optimal control law is given by (4) which yields the optimal routing in each node.
    
    We can easily verify that Assumption~\ref{as:B} and Assumption~\ref{as:A} are satisfied by the above choice of $A$, $B$ and $E$. The shortest path problem assigns cost only to traversal of the edges, while the formulation~(\ref{eq:optprob}) also associates a cost~$s_i$ with each state~$x_i$ and allows for the possibility of taking no action ($u = 0$). However, since the dynamics are not asymptotically stable~($A = I$), the optimal policy will never be to remain in a node $i$ as long as the cost $s_i$ is strictly positive. This means that the cost $s_i$ will be incurred once per edge traversal in the optimal policy. By reallocating some of the cost from the edges to all nodes except the endpoint $x_4$ (i.e., from $r$ to $s$), we guarantee that the optimal policy will correspond to the shortest path. The weights $s_i$ can be arbitrarily small. The weights \eqref{eq:weights} result in the same cost as in the original formulation for all paths that do not stay more than one time step in each node, i.e. all paths that are possible candidates for the optimal solution. 
\end{example}

\section{Distributed implementation}
Since the dimension of the state and input spaces grow in problems like \eqref{eq:optprob}, it is desirable that the complexity of the solution method scales favorably. One method for finding the optimal control of \eqref{eq:optprob} is to solve the linear program given in Theorem \ref{th:main}, which has a number of parameters that scales linearly with the dimension of the state and is thus computationally tractable even for high-dimensional problems. In many applications, however, global communication and access to information is not feasible. An alternative scalable method is to distribute the the solution of the problem, requiring only local communication. Below we present an algorithm for distributed and asynchronous calculation of the optimal local cost function for a set of agents with limited knowledge of the global state.

Consider the dynamics (\ref{eq:optprob}) when each state $x_i$ belongs to a separate agent $i$, so that $M = n$. The optimal cost function $J^*(x) = p^\top x$ with $p$ fulfilling (\ref{eq:p}) can be found by synchronous distributed value iteration starting from the ansatz $p_0 = 0$. Agent $i$ potentially affects other agents through either the dynamics~$A$, the constraints expressed in $E$ or the inputs governed by the structure of $B$. Here, let $W^{(i)}$ denote the local incidence matrix of node $i$ given by $W^{(i)} = B_iB_i^\top$. We define a set of neighbors of $i$ for each of these interactions:
\begin{align*}
    &\mathcal{N}^{A}_i := \{j: A_{ji}\neq 0, j\neq i\}\\
    &\mathcal{N}^E_i = \{j: E_{ji}\neq 0, j\neq i\}\\
    &\mathcal{N}_i^{B} = \{j: W^{(i)}_{ij}\neq 0, j\neq i\}.
\end{align*}
As a consequence of Assumption \ref{as:A} we have $\mathcal{N}_i^E \subseteq \mathcal{N}_i^A$ with equality for the choice $E = A$. Note that these sets are defined to include only the out-neighbors of $i$. Next, let us index the columns of $A$ as $A = \begin{bmatrix} A_1 ~\cdots ~A_n \end{bmatrix}$ and let each agent $i$ have local estimates $\hat{p}^{(i)},\hat{q}^{(i)}\in\mathbb{R}$. Further, let~$\hat{p}=\begin{bmatrix} \hat{p}^{(1)} ~\cdots ~\hat{p}^{(n)} \end{bmatrix}^\top$ denote the vector containing all agents' estimates of the optimal local value function. We have shown in the proof of Theorem~\ref{th:main} that~(\ref{eq:p}) has a unique solution~$p$ in~$\mathbb{R}^n_+$. Further, $\mathbb{R}_+^n$ is an invariant set of the operator given by iteration of the Bellman equation~(\ref{eq:bellman}). As a consequence, if a $p$ satisfying (\ref{eq:p}) exists, then value iteration can be performed in a distributed fashion following Algorithm \ref{alg:1} with all local estimates converging to the optimal cost function \cite[Proposition 2.6.2]{bertsekas22abstract}, under the additional assumption that all agents update their estimate an infinite number of times as~$t\to\infty$. Note that steps 5 and 6 can (with slight abuse of notation) be evaluated despite agent $i$ not having access to the full vector $\hat{p}$, since the unknown terms in both cases are multiplied by zero. As stated in Algorithm 1, each agent stores only the local estimates $\hat{p}^{(i)}$ and $\hat{q}^{(i)}$ which will converge to the fixed point corresponding to the solution of (\ref{eq:p}). In order to also extract the optimal policy, agents would additionally need to store the index of the element that minimizes the expression in step 5.

%~$\hat{p}^{(j)}$ for all $j\in \mathcal{N}_i^{B}$.

\begin{algorithm}[h!]
\caption{Asynchronous value iteration}\label{alg:1}
\begin{algorithmic}[1]
    \State local estimates $\hat{p}^{(i)}, \hat{q}^{(i)} \gets 0$ for $i$ = $1,\dots,n$
    \While{$true$}
    \State sample agent $i\in \{1,...,n\}$
    \State receive $\hat{p}^{(j)}$ from each neighbor $j \in \mathcal{N}_i^A \cup \mathcal{N}_i^B$
    \State receive $\hat{q}^{(j)}$ from each neighbor $j \in \mathcal{N}_i^E$
    \State $\hat{q}^{(i)} \gets \min\{r_i + B_i^\top \hat{p}, 0\}$
    \State $\hat{p}^{(i)} \gets s_i + A_i^\top \hat{p} + \hat{q}^{(i)}E_{ii} + \sum\limits_{j \in \mathcal{N}_i^E} \hat{q}^{(j)}E_{ji}$
    %\State broadcast $\hat{p}^{(i)}, \hat{q}^{(i)}$ to the in-neighbors of $i$, i.e. all $j$ for which $i\in\mathcal{N}_j^A \cup \mathcal{N}_j^B$ and $\mathcal{N}_j^E$ respectively
    \EndWhile
\end{algorithmic}
\end{algorithm}

% Alternative:
% \david{
% Alternatively Algorithm 2
% \begin{algorithm}[h!]
% \caption{Asynchronous value iteration}\label{alg:1}
% \begin{algorithmic}[1]
%     \State local estimates $\hat{p}^{(i)}, \hat{q}^{(i)} \gets 0$ for $i$ = $1,\dots,n$
%     \While{$true$ each agent $i$ in parallell}
%     \State receive $\hat{p}^{(j)}$ from each neighbor $j \in \mathcal{N}_i^A \cup \mathcal{N}_i^B$
%     \State receive $\hat{q}^{(j)}$ from each neighbor $j \in \mathcal{N}_i^E$
%     \State $\hat{q}^{(i)} \gets \min\{r_i + B_i^\top \hat{p}, 0\}$
%     \State $\hat{p}^{(i)} \gets s_i + A_i^\top \hat{p} + \hat{q}^{(i)}E_{ii} + \sum\limits_{j \in \mathcal{N}_i^E} \hat{q}^{(j)}E_{ji}$
%     \State broadcast $\hat{p}^{(i)}, \hat{q}^{(i)}$ to the in-neighbors of $i$, i.e. all $j$ for which $i\in\mathcal{N}_j^A \cup \mathcal{N}_j^B$ and $\mathcal{N}_j^E$ respectively
%     \EndWhile
% \end{algorithmic}
% \end{algorithm}
% }

To perform Algorithm \ref{alg:1}, agent $i$ needs to store and access, in addition to the sets of in- and out-neighbors, only the following: $A_{ji}$ for $j \in \mathcal{N}_i^A$, $E_{ji}$~for~$j~\in~\mathcal{N}_i^E$ and~$B_i$. From the perspective of privacy we note that agent $i$ requires only the quantities $\hat{p}_j$ and $\hat{q}_j~=~\min\{r_j~+~B_j^\top \hat{p}, 0\}$ from its neighbors, which do not contain information about which action attains the minimum. This means that agents are not required to reveal their optimal action in order to perform the algorithm.

\begin{remark}
    The algorithm as stated here does not include a terminating condition, which is realistic in settings where the optimal controller is calculated online and system parameters are potentially subject to change. Various methods of termination detection for distributed and asynchronous algorithms exist dating back to the seminal work \cite{huang89termination}, but applying these to the results presented here in a computationally efficient manner is a non-trivial problem in its own right and beyond the scope of this paper.
\end{remark}

\section{Example: linear flow dynamics}

We now present a more intricate, yet physically motivated, example to demonstrate some of the nuances that can be captured by the presented framework beyond the simpler shortest path problem. In the following example an optimal controller is derived for a system where endpoints are not specified a priori and the routing problem is complicated by the underlying dynamics of the network.

\begin{example}\label{ex:cool}
Consider a simplified version of a liquid cooling system serving some type of large-scale industrial plant. The system is represented by a directed graph $\mathcal{G}(\mathcal{V,E})$ where the set of links $\mathcal{E}$ models travel across the pipes that transport the liquid coolant between the nodes $\mathcal{V}$. Each node $i\in\mathcal{V}$ corresponds to a different location in the plant and has an associated state $x_i$ that specifies the total heat contained there. The dynamics are given by (\ref{eq:optprob}) with $n = |\mathcal{V}|$ and $m = |\mathcal{E}|$. Let $A_{ii}$ model the loss of heat at location $i$ due to both dissipation and diffusion to other states, while the off-diagonal elements~$A_{ij}, j\ne i$ model the heat diffusing from state $i$ to other states $j$. For all states we have $\sum_j A_{ji} \le 1$, where equality means that no significant dissipation occurs in state~$i$. 

\begin{figure}[t!]
    \centering
    \includegraphics[width=.9\linewidth]{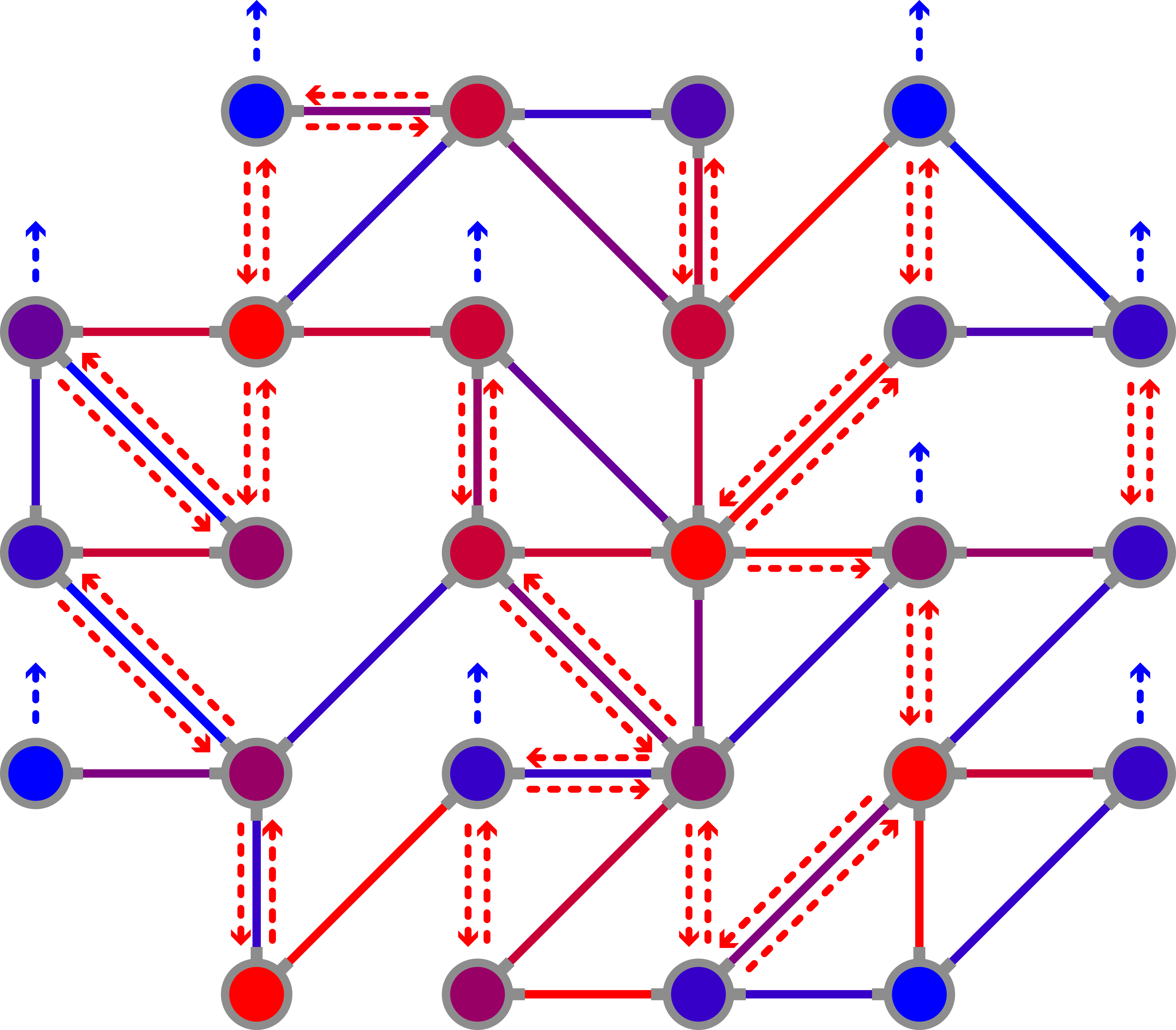}
    \caption{Graph illustrating the cooling system in the example. The coloring of nodes and edges represent the cost vectors $s$ and $r$, ranging from bright blue, corresponding to a cost of 0.2 per unit of heat transported, to bright red, corresponding to a cost of 1. Intermediate hues represent costs on that interval. Dashed red lines signify diffusion of heat between states and dashed blue lines represent dissipation.}
    \label{fig:coolnet}
    \vspace*{-2mm}
\end{figure}

\begin{figure}[t!]
    \vspace*{5mm}
    \centering
    \includegraphics[width=.9\linewidth]{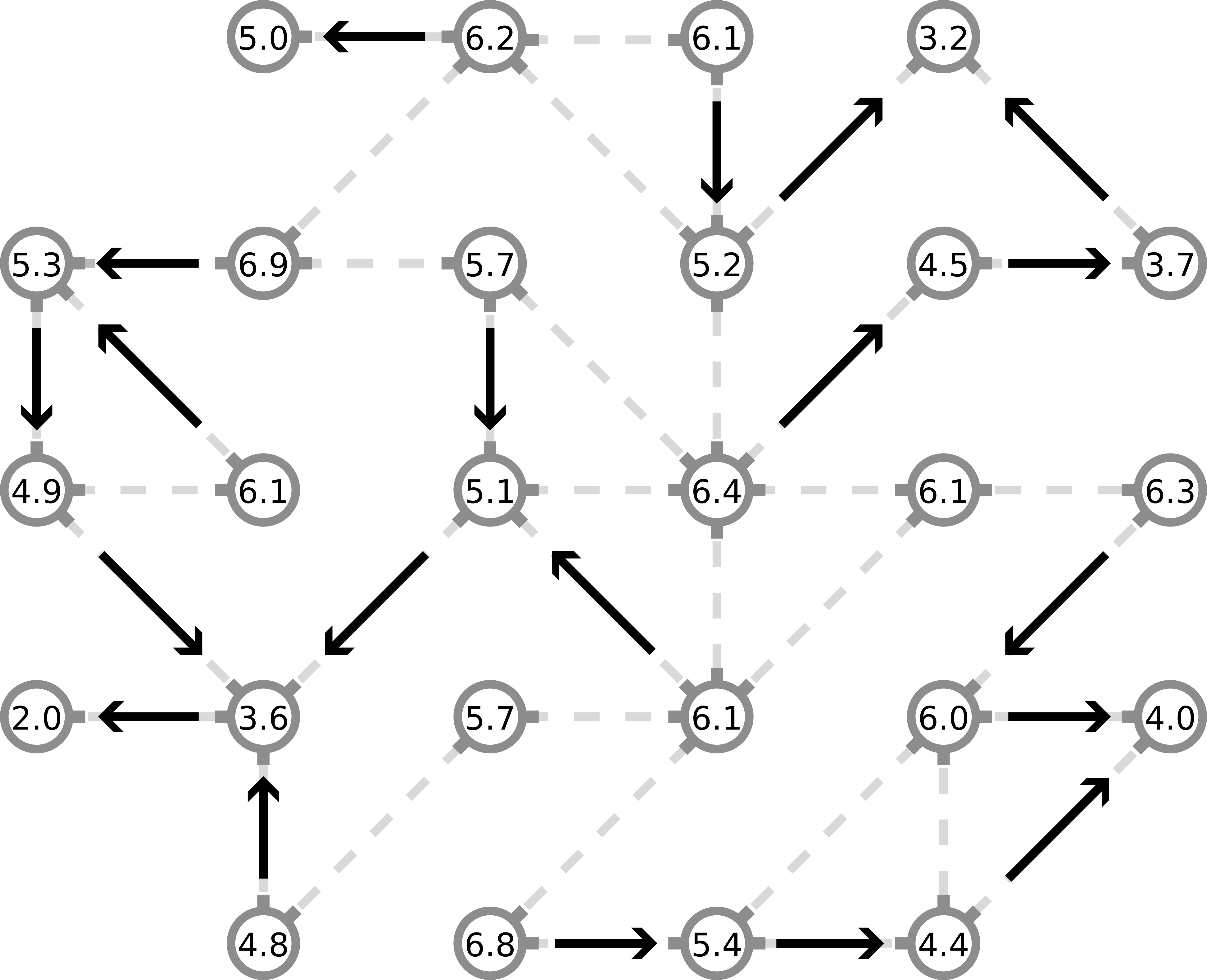}
    \caption{The optimal cost function $J^*(x) = p^\top x$ given by Theorem \ref{th:main} with each node showing the corresponding element of $p$. Dashed lines indicate the possible inputs (in both directions) and solid arrows show the direction of the optimal stationary control policy.}
    \label{fig:optnet}
    \vspace*{-2mm}
\end{figure}

In each node $i$, the heat can be transferred by selecting the input vector $u_i$, dictating how much heat is transported across each of the links~$e_{ji}$ leading from $i$. The amount transferred is limited by the current amount present in state $x_i$, giving the constraint $\mathbf{1}^\top u_i \le E_i^\top x$ where $E_{ii} = A_{ii}$ and $E_{ij} = 0$ for~$i\ne j$. The matrix $B$ has one column corresponding to each directed link $e_{ij}$ (from $j$ to $i$), in which the $j$th element is~-1 and the $i$th is a number in the interval $\left[0.95,1\right)$ modeling the heat loss over the pipe. The weights $r$ and $s$ signify the relative importance of cooling the different sections of the plant. Since each pipe is represented by two links $e_{ij}$ and $e_{ji}$, the weights and heat dissipation for the two directions are set to be identical.

In Figure \ref{fig:coolnet} the dynamics of the system are illustrated. The coloring of the nodes and links represent the relative costs $r$ and $s$ of transporting and storing one unit of heat respectively. The stationary control law found by applying Theorem~\ref{th:main} is illustrated in Figure \ref{fig:optnet} by indicating the optimal routing of heat across the graph. In those nodes $i$ that have no arrow originating in them, the optimal control is~$u_i = \mathbf{0}$, letting the heat diffuse and dissipate passively. As the dissipation only occurs in certain nodes, all heat must necessarily be transported to these either through passive diffusion or active heat transfer in order to obtain a finite value of the optimization problem~\eqref{eq:optprob}. Figure \ref{fig:cent} shows the elements of $p_k$ when applying value iteration (fixed point iteration of \eqref{eq:p}) to find the fixed point of \eqref{eq:p} for the above problem, with initialization $p_0 = \mathbf{0}$. In Figure \ref{fig:dist} the evolution of the local cost estimates $\hat{p}^{(i)}$ when performing Algorithm \ref{alg:1} for the same problem are displayed. Both methods converge to the fixed point corresponding to the optimal cost function in Figure \ref{fig:optnet} and the rate of convergence differs approximately by a factor~$n=26$, since only one agent updates their state in each iteration of Algorithm \ref{alg:1}. In other words, the figures show both methods converging after approximately the same number of agent updates.

\begin{figure}[t!]
    \centering
    \begin{tikzpicture}[trim axis left,trim axis right,baseline]

\begin{axis}[
width=.9\columnwidth,
height=6cm,
xlabel={Iterations},
ylabel={Local cost function},
xmin=0,
xmax=46,
ymin=0,
ymax=7.5,
xtick={10,20,...,70},
%extra x ticks={1},
ytick={2,4,6},
%extra y ticks={4,5},
%extra y tick labels={$\mu=4$,$\mu+\frac{\sigma^2}{\mu}$},
%extra y tick style={y tick label style={right, xshift=0.25em}},
%ytick pos=right,
%extra y tick style={grid=major,major grid style={red,line width=0.25pt},tick pos=right, ticklabel pos=right}
cycle list name=exotic,
]

\foreach \x in {1,...,26} {

    \addplot +[mark=none,solid,smooth,opacity=.7,line width=0.25pt] table[x index=0, y index=\x,col sep=comma] {pbank.csv};

}

\end{axis}

\end{tikzpicture}
    \caption{Elements of the cost function iterate $p_k$ evolving under fixed point iteration of \eqref{eq:p}, starting from $p_0 = \mathbf{0}$. Each element corresponds to a node in Example \ref{ex:cool}.}
    \label{fig:cent}
\end{figure}
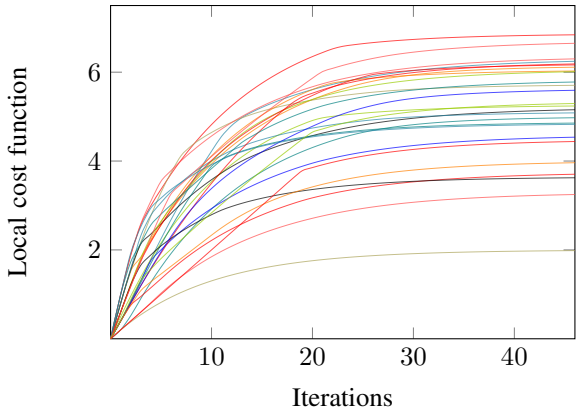

\begin{figure}[t!]
    \centering
    \begin{tikzpicture}[trim axis left,trim axis right,baseline]

\begin{axis}[
width=.9\columnwidth,
height=6cm,
xlabel={Iterations},
ylabel={Local cost function},
xmin=0,
xmax=1200,
ymin=0,
ymax=7.5,
xtick={250,500,...,1000},
%extra x ticks={1},
ytick={2,4,6},
%extra y ticks={4,5},
%extra y tick labels={$\mu=4$,$\mu+\frac{\sigma^2}{\mu}$},
%extra y tick style={y tick label style={right, xshift=0.25em}},
%ytick pos=right,
%extra y tick style={grid=major,major grid style={red,line width=0.25pt},tick pos=right, ticklabel pos=right}
cycle list name=exotic,
]

\foreach \x in {1,...,26} {

    \addplot +[mark=none,solid,smooth,opacity=.7,line width=0.25pt] table[x index=0, y index=\x,col sep=comma] {p2bank.csv};

}

\end{axis}

\end{tikzpicture}
    \caption{Local cost function estimates $\hat{p}^{(i)}$ when performing distributed value iteration according to Algorithm \ref{alg:1} for the nodes in Example \ref{ex:cool}. Note that the apparent difference in convergence rate compared to Figure \ref{fig:cent} is a result of only updating one node in each iteration of Algorithm \ref{alg:1}. Scaling the axis of iterations by a factor $n = 26$ as shown here illustrates that the rate of convergence for the two methods is similar.}
    \label{fig:dist}
\end{figure}
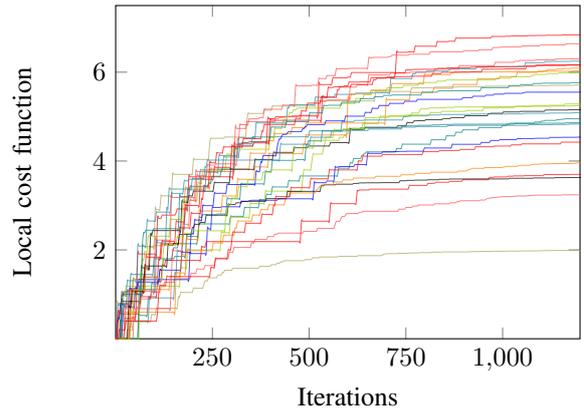

\end{example}

\section{Conclusion}

We have presented a framework for calculating optimal control of positive network dynamics in a scalable way. The resulting controller is a state feedback law~$u(t)~=~Kx(t)$ where $K$ is sparse given some inherent sparsity of the dynamics. Furthermore, we have given a distributed and asynchronous algorithm for finding the optimal cost function and corresponding controller. This framework admits a wide range of network problems, both well and lesser studied. This indicates a possible direction for further exploration in adapting tools developed on shortest path problems and Markov decision processes for use in optimal control.

\bibliographystyle{IEEEtran}
\bibliography{references}

\end{document}